\documentclass[12pt]{amsart}
\usepackage{amssymb,amscd}
\usepackage{verbatim}
\usepackage{a4wide}

\begin{document}
\newtheorem{cor}{Corollary}
\newtheorem{theorem}[cor]{Theorem}
\newtheorem{prop}[cor]{Proposition}
\newtheorem{lemma}[cor]{Lemma}
\theoremstyle{definition}
\newtheorem{defi}[cor]{Definition}
\theoremstyle{remark}
\newtheorem{remark}[cor]{Remark}
\newtheorem{example}[cor]{Example}

\newcommand{\bX}{{\partial X}}
\newcommand{\cC}{\mathcal{C}}
\newcommand{\cN}{\mathcal{N}}
\newcommand{\ch}{\mathrm{ch}}
\newcommand{\cun}{\cC^{\infty}}
\newcommand{\dom}{\mathrm{Dom}}
\newcommand{\Ind}{\mathrm{index}}
\newcommand{\ptX}{{}^{\Phi}TX}
\newcommand{\rpq}{\overline{\rz}_{+}}
\newcommand{\qz}{\mathbb{Q}}
\newcommand{\rz}{\mathbb{R}}
\newcommand{\oX}{\overline{X}}
\newcommand{\vol}{\mathrm{vol}}
\newcommand{\Xc}{X^\circ}

\title{Fibered cusp versus $d$- index theory}
\author{Sergiu Moroianu}
\thanks{Partially supported by the CERES contract 4-187/2004 
and by a CNCSIS contract (2006)}
\address{Institutul de Matematic\u{a} al Academiei Rom\^{a}ne\\ 
P.O. Box 1-764\\RO-014700 
Bucha\-rest, Romania} 
\email{moroianu@alum.mit.edu} 
\date{\today}
\begin{abstract}
We prove that the indices of fibered-cusp and $d$-Dirac operators
on a spin manifold with fibered boundary coincide if the associated
family of Dirac operators on the fibers of the boundary is invertible. 
This answers a question raised by Piazza. Under this invertibility assumption,
our method yields an index formula for the Dirac operator of
horn-cone and of fibered horn metrics.

\end{abstract}
\maketitle

Let $\oX$ be a compact manifold whose boundary is the total
space of a locally trivial fiber bundle
$\varphi:{\bX}\rightarrow Y$ of closed manifolds. 
Let $x:\oX\rightarrow\rpq$ be a defining function for
${\bX}$ and denote by $X$ the interior of $\oX$. 

The \emph{fibered cusp} (or $\Phi$-) tangent bundle $\ptX$
is a smooth vector bundle on $\oX$ defined in terms of the above data 
by its global sections:
\[\begin{split}
\cun(\oX,\ptX):=&\{V\in\cun(\oX,T\oX); V_{|\bX} \text{ is tangent to the 
fibers of $\varphi$,}\\ & \langle dx,V\rangle\in x^2\cun(\oX)\}.
\end{split}\]
When restricted to $X$, the $\Phi$ tangent bundle is canonically isomorphic 
to the usual tangent bundle $TX$. By definition, a \emph{fibered cusp} metric 
$g_\Phi$ is the restriction to $X$ of a Euclidean metric in the bundle $\ptX$, 
smooth down to the boundary of $\oX$. Let 
\[g_d:=x^2g_\Phi\]
be the conformally equivalent $d$-metric.
Such metrics appear naturally in a variety of geometric situations. 
\begin{example}
Let $(X,g)$ be a complete hyperbolic manifold of finite volume. Then outside 
a convex set, $X$  is isometric to the disjoint union of a finite number of 
``cusps'', i.e., cylinders $[0,\infty)\times M$ with metric
\[dt^2+e^{-2t}h_M,\]
where $h_M$ is flat. Compactify $X$ by setting 
$\oX:=X\sqcup (\{\infty\}\times M$). This space becomes a smooth manifold with 
boundary if we impose that $e^{-t}$ be a boundary-defining function.
By the change of variables $x=e^{-t}$, the metric becomes a 
$d$-metric for the trivial boundary fibration $M\to\{pt\}$. More generally, a  
locally symmetric space $X$ with $\qz$-rank $1$ cusps is diffeomorphic, 
outside a compact set,
to $[R,\infty)\times M$ where $M$ is the total space of a fibration 
$\phi:M\to Y$ with a canonical connection; moreover, $X$ has a 
natural Riemannian metric which near infinity takes the form
\[dt^2+\varphi^* g_Y+e^{-2t}g_Z\]
where $g_Y, g_Z$ are a metric on the base, respectively a family of metrics
on the fibers. 
\end{example}

\begin{example}
The standard metric on $X=\rz^n$ is an example of a $\Phi$ metric for $\oX$
the radial compactification, $x:=1/r$ and for the identity fibration 
at infinity $1:S^{n-1}\to S^{n-1}$ ($\Phi$-metrics
for the identity fibration are also called \emph{scattering metrics}). 
\end{example}
More generally, several examples of complete Ricci-flat metrics are of 
$\Phi$ type.

To understand better these metrics, choose a product decomposition of $X$ near 
the boundary of $\oX$ and fix a connection in the fibration $\varphi$.

\begin{defi}
A \emph{product} $\Phi$-metric (with respect to the above choices) is
a $\Phi$-metric which near infinity takes the form
\[\frac{dx^2}{x^4}+\frac{\varphi^*{g^Y}}{x^2}+h\]
where $g^Y$ is a metric on $Y$ and $h$ a family of 
metrics on the fibers, both independent of $x$. 
\end{defi}
Thus locally symmetric spaces with $\qz$-rank $1$ cusps have product 
$\Phi$-metrics.

It is straightforward to decide on the completeness of metrics conformal to
a $\Phi$-metric.
\begin{lemma}\label{complete}
Let $g$ be a $\Phi$-metric. Then for $p\in\rz$,
the metric $x^{2p}g$ is complete if and only if $p\leq 1$. 
\end{lemma}
\begin{proof}
By compactness of $\oX$, there exists a product $\Phi$-metric $g'$ on $\ptX$
and a constant $C>0$ such that 
\[C^{-1}g'\leq g\leq Cg'.\]
Hence, $x^{2p}g$ is complete if and only if $x^{2p}g'$ is. It is evident that
a metric which outside a compact set has the form
\[x^{2p-4}dx^2+h(x)\]
on $(0,\epsilon)\times M$ (where infinity corresponds to $x=0$), is complete 
if and only if the length of the segment $(0,\epsilon)\times\{m\}$ is infinite 
for all $m\in M$, or in other words if $\int_0^\epsilon x^{p-2}dx=\infty$, 
which is equivalent to $p\leq 1$.
\end{proof}

\begin{defi}
An \emph{exact} $\Phi$-metric $g_\Phi$ on $X$ is a $\Phi$-metric which differs 
from a product $\Phi$-metric by a tensor in $x \cun(\oX,S^2(\ptX))$, i.e., 
by a symmetric bilinear form on the bundle $\ptX$, smooth down to $x=0$ 
and vanishing at $x=0$. 
\end{defi}
This agrees with the definition from \cite{boris}. The name 'exact' is used 
by analogy with Melrose's exact $b$-metrics. In the rest of the paper
we will work with an exact $\Phi$-metric. These metrics were also
considered by Leichtnam, Mazzeo and Piazza \cite{lmp}, under the additional 
hypothesis that the product decomposition $[c,\infty)\times M$ is orthogonal 
\emph{near} infinity. Note that what we call $d$- respectively fibered cusp 
(or $\Phi$) metrics is called fibered cusp, respectively fibered boundary 
(or $\Phi$) metrics in \cite{lmp}. We keep our terminology for historical 
reasons.

\section{The main Theorem}

We assume that $X$, $\bX$ and the fibers of $\varphi$ have fixed (and 
compatible) spin structures. Let $E\to\oX$ be a Hermitian vector bundle 
with connection (smooth down to $x=0$). Since both $g_\Phi$ and $g_d$ 
are complete by Lemma \ref{complete}, their associated Dirac operators on $X$ 
(twisted by $E$) are essentially self-adjoint in $L^2$. 

The $L^2$ index of the Dirac operator $D^d$ was computed by Vaillant 
\cite{boris} by making extensive use of $\Phi$ operators. He assumes that 
the dimension of the kernel of the family of Dirac operators on the 
fibers of $\varphi$ is constant. In this note we observe that in the 
fully elliptic case (i.e., when the above family of Dirac operators is 
invertible) the indices of $D^d$ and $D^\Phi$ are the same. The possibility
of such a result was conjectured by Paolo Piazza in a private communication.
In general, there exists a link between the kernels, which implies
that under Vaillant's hypothesis, the index of $D^\Phi$ is finite. 

Our proof works for a more general metric $x^{2p}g_\Phi$ for some
$p\geq 0$. We denote the associated Dirac operator by $D^p$.
The metric $g_d$ is obtained by setting $p=1$, 
but interesting geometries are also obtained for other values of $p$. Using
Vaillant's work, we give at the end of the paper 
an index formula for a manifold with various such ends.

\begin{theorem}\label{th1}
Let $g_\Phi$ be an exact $\Phi$-metric on $X$, such that the 
twisted Dirac operators on the fibers of $\bX$ with respect to $E$, 
the induced metric on the fibers of the boundary 
and the induced spin structure are invertible.
Then the map of multiplication by $x^{\frac{p(n-1)}{2}}$ defines an 
isomorphism between the $L^2$-kernels of the operators $D^p$ and $D^\Phi$.
In particular the $L^2$ indices (of the chiral operators) coincide.
\end{theorem}
\begin{proof}
The spinor bundles on $X$ with respect to $x^{2p}g_\Phi$, respectively $g_\Phi$ are 
\emph{the same} (i.e., canonically identified) 
and inherit the same induced metric. We denote this unique
spinor bundle by $\Sigma$. It extends naturally to a smooth vector 
bundle over $\oX$. In the rest of the proof, we suppress the coefficient bundle
$\Sigma\otimes E$ from the notation.

\begin{lemma}\label{isom}
The unbounded operator $D^p$ acting in $L^2(X,x^{2p}g_\Phi)$ with initial domain 
$\cun_c(X)$ is unitarily equivalent to $x^{-p/2}D^\Phi x^{-p/2}$
acting in $L^2(X,g_\Phi)$ with domain $\cun_c(X)$.
\end{lemma}
\begin{proof}
These are the full (symmetric) Dirac operators.
The Dirac operators are linked by the formula of Hitchin (also used by Vaillant
\cite{boris})
\begin{equation}\label{ccdo}
D^p=x^{-p(n+1)/2} D^\Phi x^{p(n-1)/2}.
\end{equation}
The two volume forms on $X$ are related by 
\[\vol(x^{2p}g_\Phi)=x^{np}\vol (g_\Phi).\]
Thus the map
\begin{align*}\cun_c(X)&\to \cun_c(X),& \phi&\mapsto 
x^{-\frac{n}{2}}\phi
\end{align*}
is an isometry with respect to the $L^2$ inner products. The conjugation 
of $D^p$ via this isometry is illustrated by the diagram 
(of unbounded operators with domain $\cun_c(X)$):
\begin{equation}\label{diag}
\begin{CD}
L^2(X,x^{2p}g_\Phi) @>{D^p}>>
&L^2(X,x^{2p}g_\Phi)\\
@V{x^\frac{pn}{2}\cdot}VV&@V{x^\frac{pn}{2}\cdot}VV\\
L^2(X,g_\Phi)@>{x^\frac{pn}{2}D^px^{-\frac{pn}{2}}}>> 
& L^2(X,g_\Phi)
\end{CD}
\end{equation}
Using \eqref{ccdo}, we see that $x^\frac{pn}{2}D^px^{-\frac{pn}{2}}
=x^{-\frac{p}{2}}D^\Phi x^{-\frac{p}{2}}$.
\end{proof}

Unitarily equivalent operators have isomorphic kernels, hence the index of the 
chiral part $D^p_+$ of $D^p$ equals the index of 
$x^{-\frac{p}{2}}D^\Phi_+ x^{-\frac{p}{2}}$.

The analytic facts that we need from the general theory of $\Phi$ operators
\cite{mame99} are the existence
of the weighted $\Phi$-Sobolev spaces $x^{a}H^b_\Phi(X)$
on which $\Phi$-operators act, as well as the existence of parametrices 
inside the $\Phi$ calculus \cite[Proposition 8]{mame99}; that is, for every 
$A\in\Psi_\Phi^{a,b}$ elliptic there exists $B\in\Psi_\Phi^{-a,-b}$ with
\[BA-1\in\Psi_\Phi^{-\infty,0}.\]
By definition, $A\in\Psi_\Phi^{a,b}$ is called \emph{fully elliptic} 
if the normal operator of $x^b A$ is invertible as a family over $Y$ 
of suspended operators acting on the fibers of $\varphi$. 
If $A$ is fully elliptic, then $BA-1$ can even be 
made to belong to $\Psi_\Phi^{-\infty,-\infty}$.

\begin{lemma} \label{simfe}
The operators $D^\Phi$ and  
$x^{-\frac{p}{2}}D^\Phi x^{-\frac{p}{2}}$
are simultaneously fully-elliptic.
\end{lemma}
\begin{proof} Note that $D^\Phi\in\Psi_\Phi^{1,0}(X)$ and  
$x^{-\frac{p}{2}}D^\Phi x^{-\frac{p}{2}} \in \Psi_\Phi^{1,p} (X, \Sigma)$. 
In cusp-type calculi, it is a basic fact that commutation by a power of $x$
decreases the total $x$-order (i.e., \emph{increases} the power of $x$; the 
filtration is defined by the negative of the power of $x$ so that it is 
decreasing, like the symbol filtration). Namely, for $A\in 
\Psi_\Phi^{a,b}(X)$, we have 
\[[A,x^p]\in \Psi_\Phi^{a-1,b+p-1}(X).\]
This implies that the normal operator satisfies
\[\cN(x^p\cdot x^{-\frac{p}{2}}D^\Phi x^{-\frac{p}{2}})=\cN(D^\Phi).\]
\end{proof}

\begin{remark}
This lemma fails for $b$-operators, see \cite{melaps}.
\end{remark}

The proof of the following ``elliptic regularity'' lemma is standard.

\begin{lemma} \label{rapdec}
Let $A\in\Psi_\Phi^{a,b}(X,{\mathcal E},{\mathcal F})$ be fully elliptic. 
Then the $L^2$ solutions of $A\psi=0$ belong to the ideal $x^\infty 
C^\infty(\oX,{\mathcal E})$.
\end{lemma}
\begin{proof} Since $A$ is fully elliptic there exists 
$B\in\Psi_\Phi^{-a,-b}(X,{\mathcal F},{\mathcal E})$ inverting $A$ up to
$R\in x^\infty \Psi_\Phi^{-\infty}(X,{\mathcal E})$, i.e.,
\[BA=1+R.\]
Let $\psi\in L^2(X,{\mathcal E})$ be a distributional solution of the 
pseudo-differential equation $A\psi=0$. It follows
\[0=BA\psi=(1+R)\psi=\psi+R\psi\]
so $\psi=-R\psi$. But $R\in x^\infty \Psi_\Phi^{-\infty}(X,{\mathcal E})$  
implies $R\psi\in x^\infty 
C^\infty(\oX,{\mathcal E})$.
\end{proof}

We can now finish the proof of Theorem \ref{th1}. 
Recall that we assumed the family of Dirac operators $D^F$ on the fibers 
of $\varphi$ to be invertible. By \cite[Lemma 3.7]{boris} this implies
that the normal operator of $D^\Phi$ is invertible, i.e., that
$D^\Phi$ is fully elliptic. For the convenience of the reader, 
we include below a direct proof. Recall that an exact $\Phi$-metric is 
assumed to coincide with a product $\Phi$-metric up to first order terms 
in $x$. Since $x\cun(\oX,\ptX)$ is an ideal in the Lie algebra 
of fibered cusp vector fields, it follows that the Levi-Civita connection 
and the Dirac operator also agree, up to first order terms, 
with the corresponding objects for the product $\Phi$-metric.
For that metric, it is straightforward to see (once the definition of the
normal operator is recalled) that
\[\cN(D^\Phi)^2(\tau)=\|\tau\|^2_{g_Y}+(D^F)^2.\]
Therefore $\cN(D^\Phi)^2(\tau)$ is invertible for all $\tau\in TY$ 
if and only if the family $D^F$ is invertible.

In conclusion, the map 
\begin{align*}
\ker(x^{-\frac{p}{2}}D^\Phi x^{-\frac{p}{2}})\to & \ker D^\Phi&
f\mapsto& x^{-\frac{p}{2}} f
\end{align*}
is well-defined by the Lemma \ref{rapdec}, has an obviously well-defined 
inverse $f\mapsto x^\frac{p}{2} f$, and preserves parity, hence it defines
a graded isomorphism. We compose this isomorphism with
the Hilbert space isometry \eqref{diag}.
\end{proof}

The above proof holds more generally for an 
arbitrary $\Phi$-metric whose associated Dirac operator
$D^\Phi$ is fully elliptic. This condition seems however difficult to
check in practice
outside the exact case. In the cusp case (i.e., $Y$ is a point), 
a small improvement  was achieved in \cite{wlom}, where one could 
decide for a \emph{closed} cusp metric whether the Dirac operator 
is fully-elliptic or not.

For completeness of the exposition, we describe the domains of the 
operators discussed here. They are certain weighted $\Phi$-Sobolev spaces.

\begin{lemma}\label{lsa}
The operators $D^\Phi$ and $x^{-\frac{p}{2}}D^\Phi x^{-\frac{p}{2}}$
are essentially self-adjoint in $L^2_\Phi$.
The domain of the adjoint of $D^\Phi$ is $H^1_\Phi$.
If we assume moreover that $D^\Phi$ is fully elliptic, then the 
domain of the adjoint of $x^{-\frac{p}{2}}D^\Phi x^{-\frac{p}{2}}$ 
coincides with $x^pH^1_\Phi$. 
\end{lemma}

\begin{remark}
Since for $0\leq p\leq 1$ the metrics are complete, we know that the Dirac operators have precisely
one self-adjoint extension; the extra fact here for such $p$ 
is identifying the domain of the extension. 
\end{remark}

\begin{proof}
The domains of the closures of the operators 
$x^{-\frac{p}{2}}D^\Phi x^{-\frac{p}{2}}$ and $D^\Phi$ always contains
the Sobolev spaces $x^pH^1_\Phi$, respectively $H^1_\Phi$ (see \cite{mame99}). 

Since $D^\Phi$ is elliptic in the $\Phi$-sense, there exists 
$B\in\Psi_\Phi^{-1}$ with $R:= BD^\Phi-1\in\Psi_\Phi^{-\infty}$.
Let $\phi\in L^2_\Phi$ be in the domain of the adjoint of $D^\Phi$. Then
\[\phi=BD^\Phi\phi-R\phi.\]
Since by hypothesis $D^\Phi\phi\in L^2_\Phi$, we get 
$BD^\Phi\phi\in H^1_\Phi$. At the 
same time, $R\phi\in H^\infty_\Phi$. Thus $\phi\in H^1_\Phi$. Therefore
\[\dom(\overline{D^\Phi})\supset H^1_\Phi\supset \dom ({D^\Phi}^*).\]

If moreover we assume that $D^\Phi$ is fully elliptic, there exists 
$B\in\Psi_\Phi^{-1,-p}$ with $R:= Bx^{-\frac{p}{2}}D^\Phi
x^{-\frac{p}{2}}-1\in\Psi_\Phi^{-\infty,-\infty}$.
Let $\phi\in L^2_\Phi$ be in the domain of the adjoint of 
$x^{-\frac{p}{2}}D^\Phi x^{-\frac{p}{2}}$, so $\phi=
Bx^{-\frac{p}{2}}D^\Phi x^{-\frac{p}{2}}\phi-R\phi.$ By hypothesis
$x^{-\frac{p}{2}}D^\Phi x^{-\frac{p}{2}}\phi\in L^2_\Phi$, so 
$Bx^{-\frac{p}{2}}D^\Phi x^{-\frac{p}{2}}\phi\in x^pH^1_\Phi$.
Moreover, $R\phi\in x^\infty H^\infty_\Phi$, so in conclusion
$\phi\in x^pH^1_\Phi$. Thus
\[\dom \left( \overline{x^{-\frac{p}{2}}D^\Phi x^{-\frac{p}{2}}} \right)
\supset x^pH^1_\Phi\supset \dom \left({x^{-\frac{p}{2}}D^\Phi 
x^{-\frac{p}{2}}}\right)^*.\]
But for every symmetric operator $D$ we have $\dom(\overline{D})
\subset\dom(D^*)$.
\end{proof}

Since the embedding $x^{a}H^{b}_\Phi\hookrightarrow L^2_\Phi$ is compact
for $a,b>0$, it follows that fully elliptic $\Phi$-operators are 
Fredholm \cite{mame99}. Hence the kernels discussed here are all 
finite-dimensional.

In a recent work, Leichtnam, Mazzeo and Piazza \cite{lmp} computed the index 
of the Dirac operator $D^\Phi$ corresponding to an exact $\Phi$ metric for 
which the decomposition $[c,\infty)\times M$ is orthogonal near infinity. 
They deform the metric to a $b$-metric using 
results of Melrose and Rochon \cite{melroch}, apply Melrose's $b$-index 
formula \cite{melaps}, and use the adiabatic limit formula 
of Bismut and Cheeger \cite{jmj}. Together with Theorem \ref{th1} and a careful
analysis of the index density their argument gives a short proof of 
Vaillant's index formula \cite{boris} in the exact, orthogonal, 
fully elliptic case. Although the orthogonality hypothesis is used in the 
proof, the authors informed us in a private communication that the 
argument goes through under a weaker hypothesis of orthogonality only 
up to second order at the boundary.

\section{Horn-cone and fibered horn metrics}

The case $p=1$ corresponds to the $d$-metric, our original motivation,
and $p=0$ corresponds to the initial $\Phi$-metric. 
For $p>1$ the metric $x^{2p}g_\Phi$ is incomplete, nevertheless 
with the invertibility assumption along the fibers,
the associated Dirac operator is essentially self-adjoint, by Lemma \ref{lsa}. 
In this case, if we start with a product 
$\Phi$ metric, the change of variables $y=x^{p-1}$ gives (up to a constant)
a metric depending on the parameter $a=\frac{p}{p-1}>1$, 
with singularity at $y=0$:
\begin{equation}\label{hoco}
x^{2p}g_\Phi=dy^2+ y^2\varphi^*{g^Y}+y^{2a} h
\end{equation}
where $g_Y$ is a metric on the base of the boundary fibration, and $h$ is a 
family of metrics on the fibers. We call \eqref{hoco} a 
\emph{horn-cone metric}.
Similarly, for $0< p<1$, the change of variables $y=x^{1-p}$ gives for 
$a=\frac{p}{p-1}<0$:
\begin{equation}\label{fh}
x^{2p}g_\Phi=\frac{dy^2}{y^4}+ \frac{\varphi^*{g^Y}}{y^2}+y^{-2a} h
\end{equation}
We call this a \emph{fibered horn} metric. It is singular at $y=0$. 
The metric \eqref{hoco} is incomplete, while \eqref{fh} is complete,
expands with time in the base directions and shrinks in the fibers. 
In the result below, we could allow some weaker asymptotics for these metrics, 
such that after reverting back to the variable $x$ the resulting metric is 
smooth in $x$ and exact as a conformal $\Phi$-metric.
\begin{cor}
Let $(X,g)$ be a spin Riemannian manifold with a boundary
fibration $\varphi:\bX\to Y$, isometric outside a compact
set to a disjoint union of ends of the following types:
\begin{itemize}
\item  horn-cone type \eqref{hoco};
\item fibered horn \eqref{fh};
\item exact $d$-metric;
\item exact $\Phi$-metric.
\end{itemize}
The conformal weights $p$ may be different for each end.
Let $E\to X$ be an auxiliary vector bundle endowed with a Hermitian 
connection smooth in $x$ down to $x=0$, where $x=y^{\frac{1}{p-1}}$
for ends of type \eqref{hoco}, respectively $x=y^{\frac{1}{1-p}}$
for ends of type \eqref{fh}. 
If the family of (twisted) Dirac operators 
on the fibers on each end with respect to the induced
vertical family of metrics is invertible, then the Dirac operator 
$D^g$ on $X$ is essentially self-adjoint and Fredholm. 
If the dimension of $X$ is even, the index of $D^g$ is given by
\[\Ind(D^g_+)=\int_X\hat{A}(g)\ch(E)+\int_Y \hat{A}(g_Y)\hat{\eta}(D^h).\]
Moreover, $D^g$ has compact resolvent (thus, purely discrete spectrum)
if and only if $X$ has no $\Phi$-type end. In that case a generalized Weyl 
law holds as in \cite{wlom}.
\end{cor}
The right-hand side involves the eta-form defined in \cite{jmj}.
\begin{proof}
We compactify $X$ using the variable $x=y^{a-1}$ for each end of horn-cone 
type, and $x=y^{1-a}$ for ends of fibered horn type, thus our metric 
becomes of the form $x^{2p}g_\Phi$, where
$p$ is now a function from the set of ends into $[0,\infty)$. 
The metric $g_\Phi$ is exact for all the ends. 
Then by \cite[Chapter 2]{boris}, the Pontrjagin forms of 
$g$ are smooth down to $x=0$, so there is no need to regularize the 
local index integral.
The operator $D^g$ is fully elliptic by the invertibility assumption. 
By Lemma \ref{lsa}, it follows that it is essentially self-adjoint. 
The resolvent lives in $\Psi^{-1,-p}$ so it is compact if and only if $p$ 
takes values in $(0,\infty)$, i.e., if there is no $\Phi$-type end. 
By Theorem \ref{th1}, the index is the same as that of the associated 
conformal $d$-metric $g_d$, 
which is computed by Vaillant by the above formula, only that it involves the
index density for $g_d$, not for $g$. But of course, the 
Pontrjagin forms are conformal invariants so the formula follows. Alternately, 
instead of Vaillant's result we may use \cite{lmp} (see the discussion
on orthogonality at the end of the previous Section).

The Weyl law follows entirely as in \cite{wlom}, based on the fact that the
zeta function of $D^g$ has a meromorphic extension to the complex plane, with
a possibly double first pole.
\end{proof}

\section{The non-fully elliptic case}

Let us mention also a partial result in the non-fully elliptic case. 
Assume that the hypothesis of Vaillant \cite{boris} holds, namely that 
the kernels of the
family of Dirac operators along the fibers of the boundary form a vector bundle.
Then Vaillant \cite{boris} showed that the dimension of $\ker (D^d)$ is finite. 

\begin{prop}
Under the above hypothesis, for every $0\leq p\leq 1$, the kernel of the Dirac 
operator on $X$ associated to $x^{2p}g_\Phi$ is finite-dimensional.
\end{prop}
\begin{proof}
Let $D(p)$ denote the Dirac operator of $x^{2p}g_\Phi$.
Following the proof of Theorem \ref{th1}, we see that the map of 
multiplication by $x^{\frac{1-p}{2}}$ inside $L^2(X,g_d)$
maps the kernel of an operator in $L^2(X,g_d)$ unitarily
equivalent to $D(p)$, into the kernel of $D^d$. 
This multiplication map is injective but not surjective.
Thus $\ker D(p)$ injects into $\ker(D^d)$, whose dimension if finite
by \cite[Chapter 3]{boris}. 
\end{proof}
For $p=0$,
under Vaillant's hypothesis, the index of $D^\Phi$ is finite,
although from Melrose's Fredholm criterion we see that in the 
non-fully elliptic case, $D^\Phi$ is not Fredholm on any Sobolev space
with exponential weight $e^{\frac{a}{x}}x^{b}L^2(X,g_\Phi)$. 

{\small
\subsection*{Acknowledgments.} I am indebted to Paolo Piazza for 
proposing this question and for his valuable suggestions and comments.}

\bibliographystyle{plain}

\begin{thebibliography}{1}

\bibitem{jmj}
J.-M.~Bismut and J.~Cheeger,
{\sl $\eta$-invariants and their adiabatic limits, }
J.\ Amer.\ Math.\ Soc.\ {\bf 2} (1989), 33--70.

\bibitem{lmp}
E.~Leichtnam, R.~Mazzeo and P.~Piazza,
{\sl The index of Dirac operators on manifolds with fibered boundaries, }
to appear in Proc.\ Joint BeNeLuxFra Confer.\ Math.,
Ghent, May 20-22, 2005, Bull.\ Belg.\ Math.\ 
Soc. -- Simon Stevin.

\bibitem{mame99}
R.~R.~Mazzeo and R.~B. Melrose,
{\sl Pseudodifferential operators on manifolds with fibered boundaries, }
Asian J. Math. {\bf 2} (1998), 833--866.

\bibitem{melaps}
R.~B.~Melrose,
{\sl The Atiyah-Patodi-Singer index theorem, }
Research Notes in Mathematics {\bf 4}, 
A. K. Peters, Wellesley, MA (1993).

\bibitem{melroch}
R.~B.~Melrose and F.~Rochon,
{\sl Index in $K$-theory for families of fibred cusp operators, } preprint
math.DG/0507590.

\bibitem{wlom} 
S.~Moroianu,  
{\sl Weyl laws on open manifolds, } preprint 
math.DG/0310075.

\bibitem{boris}
B.~Vaillant,
{\sl Index- and spectral theory for manifolds with generalized
  fibered cusps, }
Dissertation, Bonner Math.\ Schriften {\bf 344} (2001), 
Rheinische Friedrich-Wilhelms-Universit\"at Bonn.

\end{thebibliography}

\end{document}